\documentclass{amsart}
\usepackage{amssymb}
\usepackage{amsfonts}
\usepackage{amsmath}
\usepackage{url}
\usepackage{color}

\usepackage{epsfig}
\usepackage{graphicx}
\usepackage{psfrag}

\newtheorem{theorem}{Theorem}

\newtheorem{proposition}[theorem]{Proposition}
\newtheorem{corollary}[theorem]{Corollary}
\newtheorem{lemma}[theorem]{Lemma}

\theoremstyle{remark}
\newtheorem{remark}[theorem]{Remark}

\newtheorem{definition}[theorem]{Definition}

\newcommand{\1}{1\hspace{-2.5pt}{\rm l}}
\newcommand{\binomial}[2]{\left 
(\!\!\begin{array}{c}#1\\#2\end{array}\!\!\right)}

\newcommand{\llabel}[1]{\label{#1}
}

\begin{document}

\title{Formal Markoff maps are positive}
\author[F. Gu\'eritaud]{Fran\c{c}ois Gu\'eritaud}
\date{May 2006}

\begin{abstract}This note defines a family of Laurent polynomials indexed in $\mathbb{P}^1\mathbb{Q}$ which generalize the Markoff numbers and relate to the character variety of the one-cusped torus. We describe which monomials appear in each polynomial and prove all the coefficients are positive integers. We also conjecture a generalization of that positivity result. \end{abstract}

\maketitle

\section{introduction}

In \cite{bowditch}, Bowditch defined Markoff maps as an appealing way of analyzing the length spectrum of the set $\mathcal{C}$ of simple closed geodesics on a hyperbolic one-cusped torus $S\simeq(\mathbb{R}^2-\mathbb{Z}^2)/\mathbb{Z}^2$. He noted that $\mathcal{C}$ stands in natural bijection with $\mathbb{P}^1\mathbb{Q}=\mathbb{Q}\cup\{\infty\}$ via the \emph{slope} function $$\sigma: \mathcal{C}~\widetilde{\longrightarrow}~\mathbb{P}^1\mathbb{Q}~,$$ and associated (bijectively) to each $c\in\mathcal{C}$ a complementary region $R_c$ of an infinite trivalent tree $\mathcal{T}$ properly embedded in the plane. This tree $\mathcal{T}$ is dual to the Farey triangulation of the hyperbolic plane $\mathbb{H}^2$ (see Section \ref{laurent} for definitions): namely, $R_c$ is the complementary region of $\mathcal{T}$ whose closure in the disc $\mathbb{H}^2\cup\mathbb{P}^1\mathbb{R}$ contains the ideal point $\sigma(c)$. If $\mathcal{R}$ is the collection of all the regions $R_c$, the Markoff map $$\Phi:\mathcal{R}\longrightarrow \mathbb{R}$$ associates to $R_c$ the trace of an element of $SL_2(\mathbb{R})$ representing $c$ (here we choose a lift of the holonomy representation $\pi_1(S)\rightarrow PSL_2(\mathbb{R})$).
The definition of $\Phi$ extends to Kleinian representations $\rho:\pi_1(S) \rightarrow SL_2(\mathbb{C})$, and Bowditch studied in particular the relationship between $\Phi$'s being proper and $\rho$'s being quasifuchsian. Markoff maps also provide new proofs and generalizations of McShane's identity \cite{bowditch, makmac}, and their intriguing analytic properties have not yet been fully explored.

\medskip
Of course, a Markoff map $\Phi$ is a very redundant object. It is in fact enough to know $\Phi(R_c)$ for three adjacent regions $R_c$ to reconstruct $\Phi$ completely. For instance, denote by $R_s$ the region $R_{\sigma^{-1}(s)}$ for $s\in \mathbb{P}^1\mathbb{Q}$, and consider \begin{equation} \llabel{xyz} \Phi(R_0)=X~~;~~\Phi(R_{\infty})=Y~~;~~\Phi(R_{-1})=Z.\end{equation} Then, every $\Phi(R_s)$ can be given by an explicit formula $f_s(X,Y,Z)$. There is in fact a non-trivial algebraic relationship between $X,Y,Z$, so many very different formulas for $f_s$ exist. In \cite{qf}, we were led to look for expressions of $f_s$ as a Laurent polynomial of degree $1$ in $X,Y,Z$:
\begin{equation}
f_s=\sum_{\alpha,\beta\in\mathbb{Z}} F_s(\alpha,\beta) \frac{X^{1+\alpha}\,
Y^{1+\beta}}{Z^{1+\alpha+\beta}}~~ \in \mathbb{Z}\left [X^{\pm
1},Y^{\pm 1},Z^{\pm 1} \right ] \llabel{coefficients}
\end{equation}
In Section \ref{laurent}, we show that such an expression exists, 
and that furthermore the integer $F_s(\alpha,\beta)$ equals $0$ unless $(\alpha,\beta)$ satisfies a natural parity condition. Our main theorem is

\begin{theorem} \llabel{main}
The Laurent polynomial $f_s$ has only positive coefficients. Moreover, all monomials in the Newton polygon of $f_s$ which satisfy the parity condition have nonzero coefficients.
\end{theorem}

\noindent (Recall that the Newton polygon of a Laurent polynomial $P=\sum a_{\nu_1\dots\nu_n}X_1^{\nu_1}\dots X_n^{\nu_n}$ in $n$ variables is the convex hull in $\mathbb{R}^n$ of the points $(\nu_1,\dots,\nu_n)\in\mathbb{Z}^n$ for which $a_{\nu_1\dots\nu_n}\neq 0$.) In fact, we describe the Newton polygon of $f_s$ completely (see (\ref{domain}) below). Some examples are shown in Figure \ref{abeilles} page \pageref{abeilles}. The numbers $f_s(1,1,1)$ are the usual Markoff numbers from Diophantine approximation theory \cite{refmarkoff}.

The positivity of the coefficients $F_s(\alpha,\beta)$ is already less than trivial when $s$ is a fairly simple rational of $\mathbb{P}^1\mathbb{Q}$, say an integer (that case was used in Section 7 of \cite{qf}, to establish a certain convergence property in the Teichm\"uller space of the cusped torus). In general, this author wonders about a possible interpretation (geometric, algebraic or combinatorial) of these positive numbers $F_s(\alpha,\beta)$.

\section{The functions $f_s$ are Laurent polynomials} \llabel{laurent}

Let $S=(\mathbb{R}^2-\mathbb{Z}^2)/\mathbb{Z}^2$ be the one-cusped (or once-punctured) torus and $\pi:\mathbb{R}^2-\mathbb{Z}^2\rightarrow S$ the natural projection. Denote by $\mathcal{C}$ the set of isotopy classes of simple closed curves in $S$ that are not the loop around the cusp. If $p,q$ are coprime integers and $\ell$ is a line in $\mathbb{R}^2$ of slope $s=q/p$ missing $\mathbb{Z}^2$, then $\pi(\ell)$ defines an element $c$ of $\mathcal{C}$. We call $s \in \mathbb{P}^1\mathbb{Q}$ the \emph{slope} of $c$, and write $\sigma(c)=s$. It is well-known that $\sigma$ establishes a bijection $\mathcal{C} \tilde{\rightarrow} \mathbb{P}^1\mathbb{Q}$. The curve of slope $s$ is denoted by $c_s$.

Consider the hyperbolic plane $\mathbb{H}^2$ with its natural boundary $\partial\mathbb{H}^2=\mathbb{P}^1\mathbb{R}$. Whenever two curves $c,c' \in \mathcal{C}$ have (minimal) intersecion number $1$, we connect the rationals $\sigma(c)$ and $\sigma(c')$ by a line in $\mathbb{H}^2$. The result is the \emph{Farey triangulation} of $\mathbb{H}^2$ into infinitely many ideal \emph{Farey triangles}. It is well-known that the triples of vertices of Farey triangles are exactly those triples of rationals that can be written $$\left ( \frac{q_0}{p_0}~,~\frac{q_0+q_1}{p_0+p_1}~,~\frac{q_1}{p_1} \right ) ~~~~ \text{ where } ~~~~ \left|\left(\begin{array}{cc}q_0&q_1\\ p_0&p_1\end{array}\right)\right|=\pm 1$$ (we agree that $\infty=\frac{\pm 1}{0}$). Geometrically, the Farey triangulation is generated by reflecting the triangle $1\infty 0$ in its sides \emph{ad infinitum}.

Choose a point $p\in S$.
Let $\tau$ be the trace operator on $SL_2(\mathbb{R})$, and fix a representation $\rho:\pi_1(S,p)\rightarrow SL_2(\mathbb{R})$ such that if $\gamma\in\pi_1(S,p)$ is in the conjugacy class of the loop around the puncture, then $\tau\circ\rho(\gamma)=-2$ (we say that $\rho$ is \emph{type-preserving}).

\begin{proposition}\llabel{tracerelation}
The trace $\tau$ induces a function, also noted $\tau$, on $\mathcal{C}\simeq\mathbb{P}^1\mathbb{Q}$. If $s,s_0,s_1,s'$ are elements of $\mathbb{P}^1\mathbb{Q}$ such that $s_0s_1s$ and $s_0s_1s'$ are Farey triangles, then $\tau(s)$ and $\tau(s')$ are the roots of the polynomial $X^2-\tau(s_0)\tau(s_1)X+\tau(s_0)^2+\tau(s_1)^2$. 
\end{proposition}
\begin{proof}
Defining $\tau$ on $\mathcal{C}$ is straightforward, since each curve in $\mathcal{C}$ determines a conjugacy class (together with its inverse) in the image of $\rho$. We will further omit the slope bijection $\sigma:\mathcal{C}\rightarrow\mathbb{P}^1\mathbb{Q}$ and simply consider $\tau$ as defined on $\mathbb{P}^1\mathbb{Q}$.

The modular group $SL_2(\mathbb{Z})$ acts naturally on the cusped torus $S$ while preserving the isotopy class of the loop around the cusp. The induced action on $\mathcal{C}$ coincides (\emph{via} $\sigma$) with the M\"obius action on $\mathbb{P}^1\mathbb{Q}\subset \partial \mathbb{H}^2$, which extends to an action on the Farey triangulation of $\mathbb{H}^2$ that is transitive on the set of all Farey edges $s_0s_1$.

Endow the two curves $c_{s_0},c_{s_1} \in \mathcal{C}$ with orientations and arrange $c_{s_0}$ and $c_{s_1}$ in $S$ so that they intersect only at the base point $p\in S$. Then $c_{s_0}, c_{s_1}$ define elements $g_{s_0},g_{s_1}$ of $\pi_1(S,p)$.

\emph{Observation}: $[g_{s_0},g_{s_1}]$ determines a simple loop around the puncture, and therefore has trace $-2$. The curves $c_s$ and $c_{s'}$ determine the conjugacy classes of $g_{s_0}g_{s_1}$ and $g_{s_0}g_{s_1}^{-1}$ (not necessarily in that order, depending on the chosen orientations).

This observation can be checked easily when $(s_0, s_1)=(0,\infty)$ (hence $\{s,s'\}=\{1,-1\}$). The general case follows because the curves in $\mathcal{C}$ which have intersection number $1$ with $c_{s_0}$ and $c_{s_1}$ are always exactly $c_s$ and $c_{s'}$, and the $SL_2(\mathbb{Z})$-action (transitive on Farey edges $s_0s_1$) respects the intersection numbers and the loop around the cusp. 

Recall the following trace relations, valid for all $a,b\in SL_2(\mathbb{R})$:
\begin{eqnarray*} \tau(ab)+\tau(ab^{-1})&=&\tau(a)\tau(b) \\ \tau(ab)\tau(ab^{-1})&=&\tau^2(a)+\tau^2(b)-2-\tau([a,b]).
\end{eqnarray*}
Setting $a=g_{s_0}~,~b=g_{s_1}$, the Proposition follows.
\end{proof}

In the notation above, we now define $f_s:=\tau(s)$.
Dual to the Farey triangulation is an infinite $3$-valent tree in $\mathbb{H}^2$ whose complementary regions $R_s$ stand in bijection with the Farey vertices $s\in\mathbb{P}^1\mathbb{Q}$. The Markoff map $\Phi$ is therefore defined by $\Phi(R_s)=f_s$. By Proposition \ref{tracerelation}, the variables $$(X,Y,Z)=(f_0,f_{\infty},f_{-1})$$ of (\ref{xyz}) satisfy the \emph{Markoff equation}
$$X^2+Y^2+Z^2=XYZ.$$
(This equation defines the character variety, or variety of type-preserving representations.) Moreover, Proposition \ref{tracerelation} implies that if $(A,B,C,D)=(f_{s'},f_{s_0}, f_{s_1},f_{s})$ and $A,B,C$ are known (for example in terms of $X,Y,Z$), then we can always recover $D$ by either one of the formulas $$D=BC-A \hspace{30pt} \text{ or } \hspace{30pt} D=(B^2+C^2)/A.$$ In fact, these relations allow us to define $f_s$ (and therefore $\Phi$) inductively for all $s\in \mathbb{P}^1\mathbb{Q}$, in terms of $X,Y,Z$. In order to make each $f_s=\Phi(R_s)$ a homogeneous \emph{Laurent polynomial} of degree $1$ in $X,Y,Z$, we tweak the first induction relation above and use
\begin{equation} \llabel{firstdefinition}
f_s=f_{s_0}f_{s_1}\frac{X^2+Y^2+Z^2}{XYZ}-f_{s'}\end{equation}
where $s,s_0,s_1,s'$ are as in Proposition \ref{tracerelation}. For example, $f_1=\frac{X^2+Y^2}{Z}$. For all $s\in\mathbb{P}^1\mathbb{Q}$, denote by $[s]$ the unique element of $\{0,-1,\infty\}$ such that $s$ and $[s]$ project to the same point of $\mathbb{P}^1(\mathbb{Z}/2\mathbb{Z})$. In particular, $f_{[s]}$ is one of the variables $X,Y,Z$.

\begin{proposition}
If $f_s$ is defined inductively for all $s\in\mathbb{P}^1\mathbb{Q}$ using (\ref{firstdefinition}), then $f_s$ is a Laurent polynomial in $X,Y,Z$. Moreover there is a finitely supported function $F_s:\mathbb{Z}^2\rightarrow \mathbb{Z}$ such that $$f_s=\left ( \sum_{\alpha,\beta\in\mathbb{Z}} F_s(\alpha,\beta) \frac{X^{1+\alpha} \,Y^{1+\beta}}{Z^{1+\alpha+\beta}}\right )~~ \in f_{[s]} \cdot \mathbb{Z}\left [X^{\pm 2},Y^{\pm 2},Z^{\pm 2} \right ].$$ \llabel{evendegrees}
\end{proposition}
\begin{proof} From (\ref{firstdefinition}), 
by induction, $f_s$ is a Laurent polynomial. The claim on the parity of the degrees also follows by induction from (\ref{firstdefinition}), because $\{f_{[s_0]},f_{[s_1]},f_{[s]}\}=\{X,Y,Z\}=\{f_{[s_0]},f_{[s_1]},f_{[s']}\}$ holds whenever $s,s_0,s_1,s'$ are as in Proposition \ref{tracerelation}.
\end{proof}

In Section \ref{durdur} we prove Theorem \ref{main} for positive rationals $s$. The remaining cases ($s<-1$ and $-1<s<0$) will follow by a symmetry argument (see Section \ref{grozarbre}). Section \ref{generalization} exposes a generalization of our ``tweaking'' operation (\ref{firstdefinition}), and a conjecture extending Theorem \ref{main}.

\section{A family of domains and functions} \llabel{durdur}

Define $\mathcal{Q}=\mathbb{Q}^{\geq 0}\cup\{\infty\}$. Any point $s$ of $\mathcal{Q}$ can be written in a unique way $$s=\frac{q}{p} ~~~~\text{ with $p,q\in\mathbb{N}$ coprime}$$ (we agree that $\infty=\frac10$). For such $s \in
\mathcal{Q}$, define \begin{equation} \llabel{domain} J_s:=\left \{
(\alpha,\beta)\in \mathbb{Z}^2 \left | \begin{array}{l} \alpha\equiv q~;~
\beta\equiv p ~~[2] \\ \alpha\geq -q~;~\beta\geq -p \\ \alpha+\beta \leq p+q-2
\\ p\alpha+q\beta \geq 0 \end{array} \right . \right \}. \end{equation} It will
turn out that $F_s$ is supported exactly on $J_s$. Observe that $J_0=\{(0,-1)\}$
and $J_{\infty}=\{(-1,0)\}$ and $J_1=\{(-1,1);(1,-1)\}$. Further, define
\begin{itemize}
\item $Z_s=(q,p)+2\mathbb{Z}^2$ so that $J_s\subset Z_s$ ;
\item $P^s_i=(q+2i,-p)~\in Z_s$ for all $i\in\mathbb{Z}$ ;
\item $Q^s_j=(-q,p+2j)~\in Z_s$ for all $j\in\mathbb{Z}$ ;
\item $\varphi_s(\alpha,\beta)=p\alpha+q\beta$ ;
\item $\Lambda=\{(0,0);(0,2);(2,0)\}$ ; 
\item $n\Lambda=\Lambda+\dots+\Lambda=\{(2i,2j)\in 2\mathbb{N}^2|i+j\leq n\}$ 
for all $n\in\mathbb{N}$ ;
\item If $U$ is a subset of $Z_s$, then $\langle U \rangle_s$ denotes the
intersection with $Z_s$ of the convex hull of $U$ in $\mathbb{R}^2$.
\end{itemize}
\begin{figure}[h!]
\centering \input{vertices.pstex_t} \caption{The domain $J_s$. \llabel{vertices}}
\end{figure}
\begin{lemma} For all $s$ in $\mathcal{Q}$, one has $P^s_{p-1},Q^s_{q-1}\in J_s$
and $$J_s= \left \langle
\{P^s_i\,|\,0\leq i <p\}\cup\{Q^s_j\,|\,0\leq j < q\}\right \rangle_s.$$
\llabel{shapelem} \end{lemma}
\begin{proof} Having checked the two cases $s=0,\infty$ separately (one of the
families $\{P^s_i\},\{Q^s_j\}$ is then empty, so the second statement does not imply the first), assume $p,q\geq 1$ and focus on the second statement. Observe that $P^s_{p-1}, P^s_0, Q^s_0, Q^s_{q-1}$ are (in that order) the extremal points of a convex quadrilateral (or triangle, or segment, when $p=1$ and/or $q=1$), as shown in Figure \ref{vertices} (left). The sides of the quadrilateral
correspond to the four inequalities defining $J_s$, hence the result.
\end{proof}
\begin{corollary} For all $s$ in $\mathcal{Q}$ and $n$ in $\mathbb{N}$, one
has \begin{eqnarray*} J_s+n\Lambda&=&\left \langle \{P^s_i | 0\leq i
<p+n\}\cup\{Q^s_j|0\leq j < q+n\}\right \rangle_s \\ J_s+\Lambda &\supset&
[P^s_0+p\Lambda] \cup [Q^s_0+q\Lambda]. \end{eqnarray*} \llabel{shapelem2}
\end{corollary}
\begin{proof} Again, check the cases $s=0,\infty$ separately. If $p,q\geq 1$,
the first statement follows easily from Lemma \ref{shapelem} (which covers the
case $n=0$), and the second follows from the first (with $n=1$) by observing
that $P^s_0+p\Lambda$ and $Q^s_0+q\Lambda$ are the convex hulls of points of
$J_s+\Lambda$: for instance, \begin{eqnarray*} P^s_0+p\Lambda&=&\left \langle
\left \{ P^s_0;P^s_p;(q,p)\right \} \right \rangle_s \\ &=&\left \langle
\left \{ P^s_0;P^s_p; \frac{qP^s_p+pQ^s_q}{q+p}\right \} \right \rangle_s~.
\end{eqnarray*}\end{proof}

We now redefine the coefficient functions $F_s(\cdot,\cdot)$ of Proposition \ref{evendegrees} from a slightly altered point of view.
Let $\mathcal{F}$ be the $\mathbb{Z}$-module of functions $F:\mathbb{Z}^2\rightarrow \mathbb{Z}$ having finite support. We can define a convolution law on $\mathcal{F}$ by $F*G(u)=\sum_{x+y=u}F(x)G(y)$. Also, denoting by $\1_U$ the
characteristic function of a set $U$, define the following elements of
$\mathcal{F}$: $$F_s=\1_{J_s} \text{ for } s \in \{0,1,\infty\}~.$$

It is straightforward to check that the identity of Proposition \ref{evendegrees} holds for $s\in\{0,1,\infty\}$.
Finally, for $s\in \mathcal{Q}-\{0,1,\infty\}$, we shall define $F_s$ in an
inductive way. In $\mathbb{H}^2$ endowed with the Farey triangulation, consider
the line $L_s$ connecting $s$ to the midpoint $\sqrt{-1}$ of the line
$0\infty$. Denote by $s_0,s_1$ the ends of the first Farey edge encountered by
$L_s$ (closest to $s$). We call $s_0$ and $s_1$ the \emph{parents} of $s$. Up
to exchanging indices, we may assume that the parents of $s_1$ are $s_0$ and
another point $s'\in \mathcal{Q}$ (we agree that the parents of $1$ are $0$
and $\infty$). See Figure \ref{case12}. In particular, one has \begin{equation} \llabel{lollipop} \left
. \begin{array}{rcl} (p,q)&=&(p_1,q_1)+(p_0,q_0) \\
(p',q')&=&(p_1,q_1)-(p_0,q_0) \end{array} \right \} \text{ for }
\textstyle{(s,s',s_0,s_1)=(\frac{q}{p}, \frac{q'}{p'}, \frac{q_0}{p_0},
\frac{q_1}{p_1})~}.\end{equation}

\begin{definition} For each configuration as above, we set
\begin{equation} \llabel{convolez} F_s:=(F_{s_0}*F_{s_1}*\1_{\Lambda}) -
F_{s'}~~~~\text{ where }\Lambda=\{(0,0);(0,2);(2,0)\} .\end{equation}
\end{definition}

Since the dual of the Farey triangulation is a tree, this definition is easily
seen to be consistent. Clearly, $F_s$ is in $\mathcal{F}$. It is easy to check
that (\ref{convolez}) is just a reformulation of (\ref{firstdefinition}), so (\ref{convolez}) agrees with our first definition (Prop. \ref{evendegrees}) of $F_s$. The following three Lemmas (numbered \ref{fourmiliere}-\ref{extremal}-\ref{bunch}) are intended to prove that $F_s$ is supported on $J_s$ and $F_s(J_s)>0$, for all $s\in\mathcal{Q}$. The reader is invited to read their three statements first (the three proofs could be written as one vast simultaneous induction on $s$ for the simultaneous three statements).

\begin{lemma} For each configuration as above where $s\in\mathcal{Q}-\{0,1,\infty\}$, the set $J_{s'}\smallsetminus J_s$ consists of a unique (extremal) point $x_s$ of $J_{s'}$, and $J_{s_0}+J_{s_1}+\Lambda
=J_s\sqcup\{x_s\}$. \llabel{fourmiliere}\end{lemma} \emph{Remark:} if
$s\in\mathcal{Q}-\{0,\infty\}$, following Lemma \ref{shapelem}, we call
``extremal'' the points $P^s_0, P^s_{p-1}, Q^s_0, Q^s_{q-1}$ of $J_s$ (with
possible repeats). If $s\in\{0,\infty\}$, then $J_s$ is reduced to an
(extremal) point $P^s_{p-1}=Q^s_{q-1}$. \begin{proof} Let
$(\alpha,\beta)$ be an element of $J_{s'}$.  By (\ref{lollipop}) one has
$Z_{s'}=Z_s$ so $(\alpha,\beta)$ satisfies the congruence conditions of
(\ref{domain}). Still by (\ref{lollipop}), one has $p'\leq p$ and $q'\leq q$ so
the first three inequalities of (\ref{domain}) are also satisfied at
$(\alpha,\beta)$. For the fourth inequality, consider the linear form
$\varphi_s(\alpha,\beta)=p\alpha+q\beta$. Clearly, $\varphi_s(Z_s)\subset
2\mathbb{Z}$. Furthermore, observe
\begin{eqnarray*}\varphi_s(P^{s'}_i)&=&pq'-qp'+2ip\\ \varphi_s(Q^{s'}_j)
&=&qp'-pq'+2jq \\pq'-qp'&=&2(p_0q_1-p_1q_0)~=~\pm 2~~~\text{($s_0,s_1$ Farey
neighbors).}\end{eqnarray*} Thus, if $p'=0$ (resp. $q'=0$), taking for $x_s$
the only point $Q^{s'}_0$ (resp. $P^{s'}_0$) of $J_{s'}$ yields
$\varphi_s(x_s)=-2$. If $p'q'>0$, we find that exactly one point $x_s$ among
$\{P^{s'}_0,Q^{s'}_0\}$ satisfies $\varphi_s(x_s)=-2$ while $\varphi_s(x)\geq
0$ at all other extremal points $x$ of $J_{s'}$. It follows that on
$J_{s'}-\{x_s\}$ one has $\varphi_s>-2$ i.e. $\varphi_s\geq 0$. Hence the first
statement.

Let us now prove the second statement. For $(y_0,y_1,\lambda)\in J_{s_0}\times
J_{s_1} \times \Lambda$, it is again straightforward to check that
$(\alpha,\beta)=y_0+y_1+\lambda$ satisfies the congruence conditions and the
first three inequalities of (\ref{domain}). For the fourth, compute
$$\begin{array}{rcl}
\varphi_s(P^{s_0}_i)=p_1q_0-p_0q_1+2ip&&
\varphi_s(P^{s_1}_i)=p_0q_1-p_1q_0+2ip\\
\varphi_s(Q^{s_0}_j)=p_0q_1-p_1q_0+2jq&&
\varphi_s(Q^{s_1}_j)=p_1q_0-p_0q_1+2jq.\\ \end{array}$$
Again, observe that $p_0q_1-p_1q_0=\pm 1$. The same argument as above
(involving this time extremal points of $J_{s_0},J_{s_1}$ instead of $J_{s'}$) shows that $\varphi_s$ takes the value $-1$ at exactly one point
$y_0\in\{P^{s_0}_0,Q^{s_0}_0\}$ (resp. $y_1\in\{P^{s_1}_0,Q^{s_1}_0\}$) and
$\varphi_s\geq 1$ holds on $J_{s_0}-\{y_0\}$ (resp. $J_{s_1}-\{y_1\}$).
Moreover, $y_k$ belongs to $J_{s_k}$ for $k\in\{0,1\}$ (this is immediate from Lemma \ref{shapelem}, unless $p_kq_k=0$
where we need to check separately). The following
table summarizes the two possible cases for $y_0,y_1,x_s$.
\begin{equation}\begin{array}{c|c||c|c|c} \llabel{whosxs}
&p_0q_1-p_1q_0 &y_0&y_1&x_s  \\ \hline
\text{Case } 1&-1&Q^{s_0}_0&P^{s_1}_0&P^{s'}_0\\
\text{Case } 2&~~1&P^{s_0}_0&Q^{s_1}_0&Q^{s'}_0
\end{array}\end{equation}
\begin{figure}[h!] \centering
\psfrag{0}{$0$}
\psfrag{s0}{$s_0$}
\psfrag{s1}{$s_1$}
\psfrag{s}{$s$}
\psfrag{sp}{$s'$}
\psfrag{c1}{Case $2$}
\psfrag{c2}{Case $1$}
\includegraphics{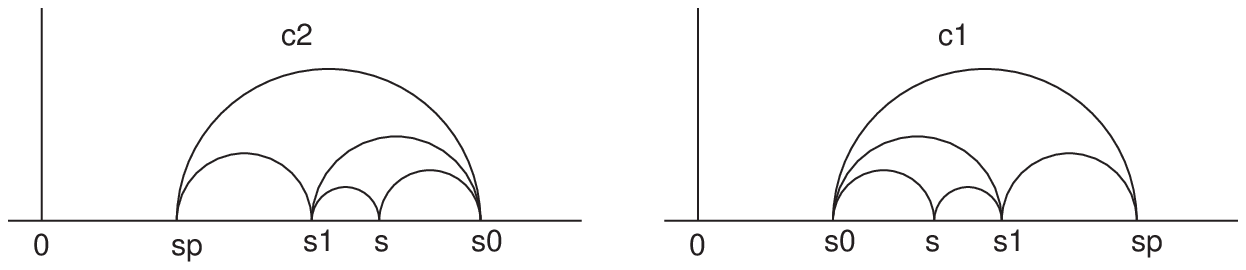}
\caption{\llabel{case12}} \end{figure}

Using Relations (\ref{lollipop}) and the definitions of $P^s_i$ and $Q^s_j$, one checks immediately that $y_0+y_1=x_s$ in both cases.
Since $\varphi_s$ is linear, $x_s$ turns out to be the only point of
$J_{s_0}+J_{s_1}+\Lambda$ where $\varphi_s<0$. This gives  one inclusion of
the equality to be proved.

For the other inclusion, $J_s \sqcup\{x_s\} \subset J_{s_0}+ J_{s_1}+ \Lambda$,
we shall restrict to Case $1$ above (Case $2$ is similar).
By Table (\ref{whosxs}), since $Q^{s_0}_0$ and $P^{s_1}_0$ belong to $J_{s_0}$ and $J_{s_1}$, one has $q_0,p_1>0$. In view of Corollary \ref{shapelem2}, it is
sufficient to prove that \begin{equation} J_s\sqcup\{x_s\}\subset (J_{s_0}+
P^{s_1}_0+ p_1\Lambda) \cup (J_{s_1}+Q^{s_0}_0+q_0\Lambda).\llabel{coverjs}
\end{equation} Still by Corollary \ref{shapelem2}, since
$P^{s_1}_0+Z_{s_0}=Z_s$, one has \begin{eqnarray*} J_{s_0}+
P^{s_1}_0+p_1\Lambda&=& \left\langle P^{s_1}_0+ \left ( \{P^{s_0}_i |0\leq
i<p_0+p_1\} \cup\{Q^{s_0}_j |0\leq j<q_0+p_1\} \right )\right\rangle_s \\ &=&
\left\langle \{P^s_i|0\leq i<p\} \cup \{Q^{s_0}_0+ P^{s_1}_0,
Q^{s_0}_{q_0+p_1-1}+ P^{s_1}_0\}\right \rangle_s \\ &=& \left \langle
\{P^s_i|0\leq i<p\}\cup \{x_s,T_0\}\right\rangle_s \\ \text{ where }
T_0&=&(q-2q_0,p+2(q_0-1)).  \end{eqnarray*} (To write the second line, we
replaced the collection of the $Q^{s_0}_j$ by its extremal terms: this is
justified because $q_0+p_1>0$). Similarly, \begin{eqnarray*}J_{s_1}+
Q^{s_0}_0+q_0\Lambda&=&\left\langle \{Q^s_j|0\leq j<q\}\cup\{x_s, T_1 \}\right
\rangle_s \\ \text{ where } T_1&=&(q+2(p_1-1),p-2p_1).\end{eqnarray*} We just
captured all the $P^s_i,Q^s_j$ which according to Lemma \ref{shapelem} define $J_s$ (Figure \ref{vertices}, right). Observe that $T_0$ (resp. $T_1$) has the same abscissa (resp. ordinate)
as $x_s=(q',-p')$. Finally, the facts that the points $Q^s_{q-1},
T_0,T_1,P^s_{p-1}$ lie in that order on the edge $E=Q^s_{q-1}P^s_{p-1}$ of
$J_s$, and that the edge $P^s_0 Q^s_0$ of $J_s$ (defined by ``$\varphi_s=0$'')
separates $x_s$ from $E$, imply (\ref{coverjs}). See the right panel of Figure \ref{vertices}. \end{proof}

\begin{lemma} \llabel{extremal} The function $F_s$ is supported on a subset of $J_s$ for all $s\in\mathcal{Q}$, and if $c$ is an extremal point of $J_s$, then
$F_s(c)=1$. \end{lemma}
\begin{proof} We prove both facts by simultaneous induction. They hold for
$s\in\{0,1,\infty\}$ so assume they hold for $s_0,s_1,s'$ and let us prove them
for $s$. By (\ref{convolez}), $F_s$ is supported on $(J_{s_0}+J_{s_1}
+\Lambda)\cup J_{s'}=J_s\cup\{x_s\}$, with $x_s$ defined as in Lemma
\ref{fourmiliere}. Recall the linear form $\varphi_s$ from the proof of Lemma
\ref{fourmiliere}: over $J_{s_0}, J_{s_1}, \Lambda$, the form $\varphi_s$
achieves its respective minima only at the extremal points $y_0, y_1, 0$;
therefore $x_s$ is realized in $J_{s_0}+J_{s_1}+\Lambda$ only as
$y_0+y_1+(0,0)$. Hence, by induction, $F_{s_0}*F_{s_1}*\1_{\Lambda}(x_s)=1$. But
$x_s$ is also an extremal point of $J_{s'}$, so (\ref{convolez}) yields
$F_s(x_s)=0$: the function $F_s$ is supported within $J_s$.

Next, observe that the extremal point $P^s_{p-1}$ of $J_s$ maximizes the first
coordinate (a similar statement is true for $J_{s_0},J_{s_1},J_{s'}$). Since
$P^s_{p-1}=P^{s_0}_{p_0-1}+P^{s_1}_{p_1-1}+(2,0)$, one has by induction
$F_{s_0}*F_{s_1}*\1_{\Lambda}(P^s_{p-1})=1$. Also, $P^s_{p-1}$ does not belong
to $J_{s'}$ because all $(\alpha,\beta)$ in $J_{s'}$ satisfy $\alpha+\beta\leq
p'+q'-2<p+q-2$. By (\ref{convolez}), we find $F_s(P^s_{p-1})=1$. Similarly,
$F_s(Q^s_{q-1})=1$. Consider one of the (at most two) remaining extremal
points of $J_s$, say $P^s_0$. Without loss of generality, one has $p\geq 2$
(otherwise, the point has already been treated as $P^s_{p-1}$). One cannot have
$\{p_0,p_1\}=\{0,p\}$ lest $|p_0q_1-p_1q_0|\geq p>1$ (recall $s_0,s_1$ are
Farey neighbors). Therefore $p_0,p_1\geq 1$. Observe that the points
$P^s_0,P^{s_0}_0,P^{s_1}_0$ are the minimizers over $J_s,J_{s_0},J_{s_1}$ of
the form $(\alpha,\beta)\mapsto \beta+\varepsilon\alpha$, for very small
$\varepsilon$. Since $P^s_0=P^{s_0}_0+P^{s_1}_0+(0,0)$, we find that
$F_{s_0}*F_{s_1}*\1_{\Lambda}(P^s_0)=1$. Finally, $P^s_0$ cannot belong to
$J_{s'}$ because of its second coordinate, $-p<-p'$. By (\ref{convolez}), this
yields $F_s(P^s_0)=1$. Similarly, $F_s(Q^s_0)=1$. \end{proof}

\begin{lemma} \llabel{bunch} For all $s\in\mathcal{Q}$ one has $F_s(J_s)\subset\mathbb{Z}^{>0}$. If $s\notin\{0,\infty\}$ then
$$\1_{J_s}\cdot\sup\left\{ \begin{array}{rcl} \1_{\{P^{s_0}_0\}}*F_{s_1} &,&
\1_{\{P^{s_1}_0\}}*F_{s_0}~~, \\ \1_{\{Q^{s_0}_0\}}*F_{s_1} &,&
\1_{\{Q^{s_1}_0\}}*F_{s_0} \end{array} \right\}\leq F_s.$$  \end{lemma}
\begin{remark} By Corollary \ref{shapelem2} and Lemma \ref{fourmiliere}, each function in the bracket is supported within $J_s\sqcup\{x_s\}$ (because e.g. $P^{s_0}_0\in J_{s_0}+\Lambda$). In other words, $\1_{J_s}$ can be replaced by $\1_{Z_s-\{x_s\}}$ without altering the strength of the statement.\end{remark}
\begin{proof} Again, both facts are proved by simultaneous induction. They
hold for $s\in\{0,1,\infty\}$; assume they hold for $s_0,s_1,s'$; let us prove
them for $s$. Recall our convention that the parents of $s_1$ are $s_0$ and
$s'$ (so in particular, $s_1\neq 0,\infty$). We saw in the course of proving
Lemma \ref{fourmiliere} that $x_s$ is either $P^{s_0}_0+Q^{s_1}_0=Q^{s'}_0$ or
$Q^{s_0}_0+P^{s_1}_0=P^{s'}_0$. On the other hand, $x_{s_1}$ is either
$P^{s_0}_0+Q^{s'}_0$ or $Q^{s_0}_0+P^{s'}_0$. In fact, using (\ref{lollipop})
and the generic characterization $\varphi_{\sigma}(x_{\sigma})=-2$, it is easy
to check that \begin{equation}\begin{array}{ccccc}
x_{s_1}=P^{s_0}_0+Q^{s'}_0&\iff& q_0p'-p_0q'=-1& \iff & x_s=Q^{s'}_0~; \\
x_{s_1}=Q^{s_0}_0+P^{s'}_0&\iff& p_0q'-q_0p'=-1& \iff &
x_s=P^{s'}_0~.\end{array} \llabel{xspq}\end{equation} Define in general
$G_s=F_s*\1_{\Lambda}$. Lemma \ref{extremal} easily yields
$G_{\sigma}(P^{\sigma}_0)=G_{\sigma}(Q^{\sigma}_0)=1$ for all
$\sigma\in\mathcal{Q}$ (this should again be checked separately for
$\sigma=0,\infty$). By Lemma \ref{fourmiliere} and the
induction hypothesis, we have $F_{s_0}*F_{s_1}*\1_{\Lambda}>0$ on $J_s$. Moreover, by (\ref{convolez}), 
\begin{eqnarray*} F_s+F_{s'}&=&F_{s_0}*F_{s_1}*\1_{\Lambda}\\&=&\sum_{\lambda\in
(J_{s_0}+\Lambda)} G_{s_0}(\lambda)\cdot \1_{\{\lambda\}}*F_{s_1} \\
&=&\left[\left(\1_{\{P^{s_0}_0\}}+\1_{\{Q^{s_0}_0\}}\right)*F_{s_1}\right ]+
\underset{\lambda \neq P^{s_0}_0, Q^{s_0}_0}{\sum_{\lambda \in
(J_{s_0}+\Lambda)}} G_{s_0}(\lambda)\cdot \1_{\{\lambda\}}*F_{s_1}~.
\end{eqnarray*} Thus, if we prove
\begin{equation}\1_{\{P^{s_0}_0\}}*F_{s_1}(x)\geq F_{s'}(x)~~;~~
\1_{\{Q^{s_0}_0\}}*F_{s_1}(x)\geq F_{s'}(x) ~~~\text{ for all }x\neq x_s~,
\llabel{shifts}\end{equation} then we will have at once $F_s>0$ on $J_s$ (because $F_s+F_{s'}\geq 2F_{s'}$ and $F_{s'}(J_{s'})>0$), and also
$F_s\geq\sup\left\{ \1_{\{P^{s_0}_0\}}*F_{s_1}, \1_{\{Q^{s_0}_0\}}*F_{s_1}
\right\}$ on $J_s$. That is half of Lemma \ref{bunch}.

Using the relation $P^{s_0}_0=-Q^{s_0}_0$ and the identities
$\1_{\{\xi\}}*\1_{\{\eta\}} =\1_{\{\xi+\eta\}}$ and
$\1_{\{\xi\}}*f(x+\xi)=f(x)$, Equation (\ref{shifts}) is equivalent to
\begin{eqnarray}F_{s_1}(y)&\geq& \1_{\{Q^{s_0}_0\}}*F_{s'}(y) ~~~\text{ if
}~~~ y\neq x_s+Q^{s_0}_0 \llabel{xsp}\\ F_{s_1}(y)&\geq&
\1_{\{P^{s_0}_0\}}*F_{s'}(y) ~~~\text{ if }~~~y\neq x_s+ P^{s_0}_0.
\llabel{xsq} \end{eqnarray} 
For $y\neq x_{s_1}$, both inequalities are already true by induction ($s_0,s'$ are the parents of $s_1$). For $y=x_{s_1}$, in view of (\ref{xspq}), two cases may arise: \begin{itemize}
\item If $x_s=P^{s'}_0$ then $x_{s_1}=x_s+Q^{s_0}_0$ so (\ref{xsp}) is true,
and (\ref{xsq}) need only be checked at $y=x_{s_1}$. One has
$F_{s_1}(x_{s_1})=0$ and $$\1_{\{P^{s_0}_0\}}*F_{s'}(x_{s_1})
=F_{s'}(x_{s_1}-P^{s_0}_0) =F_{s'}(P^{s'}_0+2Q^{s_0}_0).$$ However,
(\ref{lollipop}) yields $\varphi_{s'}(P^{s'}_0+2Q^{s_0}_0)=2(p_0q'-q_0p')=-2$;
hence, the point $(P^{s'}_0+2Q^{s_0}_0)$ does not belong to $J_{s'}$ and
$\1_{\{P^{s_0}_0\}}*F_{s'}(x_{s_1})=0$. \item Similarly, if $x_s=Q^{s'}_0$
then (\ref{xsq}) is true, and for (\ref{xsp}) one need only check
$\1_{\{Q^{s_0}_0\}}*F_{s'}(x_{s_1})=F_{s'}(Q^{s'}_0+2P^{s_0}_0)=0$ because
$\varphi_{s'}(Q^{s'}_0+2P^{s_0}_0)=-2<0$. \end{itemize}

It remains to prove $F_s\geq\sup\left\{ \1_{\{P^{s_1}_0\}}*F_{s_0},
\1_{\{Q^{s_1}_0\}}*F_{s_0} \right\}$ on $J_s$. By the lower bounds on $F_s$ we just established, it is enough to make sure \begin{equation}\llabel{s0s1}
\1_{\{P^{s_1}_0\}}*F_{s_0} \leq \1_{\{P^{s_0}_0\}}*F_{s_1} ~~;~~
\1_{\{Q^{s_1}_0\}}*F_{s_0} \leq \1_{\{Q^{s_0}_0\}}*F_{s_1} ~~ \text{ on }
Z_s-\{x_s\}.\end{equation} We focus only on the first inequality (the second
is similar). It is equivalent, by the same method as above, to:
$$\1_{\{P^{s'}_0\}}*F_{s_0}(y)\leq F_{s_1}(y) ~~~\text{ if }y\neq
x_s+Q^{s_0}_0 $$ (we used $P^{s'}_0=P^{s_1}_0+Q^{s_0}_0$, a consequence of
(\ref{lollipop})). But that inequality is true (by induction) as long as $y\neq x_{s_1}$. Again, in view of (\ref{xspq}), two cases may arise at $y=x_{s_1}$: \begin{itemize}
\item If $x_s=P^{s'}_0$ then $x_{s_1}=x_s+Q^{s_0}_0$ and there is nothing to
do; \item If $x_s=Q^{s'}_0$ we only need check the inequality above at
$y=x_{s_1}$. On one hand, $F_{s_1}(x_{s_1})=0$; on the other,
$$\1_{\{P^{s'}_0\}}*F_{s_0}(x_{s_1})=F_{s_0}(x_{s_1}-P^{s'}_0)
=F_{s_0}(P^{s_0}_0+2Q^{s'}_0)$$ but, by (\ref{lollipop}),
$\varphi_{s_0}(P^{s_0}_0+2Q^{s'}_0)=2(q_0p'-p_0q')=-2<0$ so the point
$(P^{s_0}_0+2Q^{s'}_0)$ does not belong to $J_{s_0}$ and
$\1_{\{P^{s'}_0\}}*F_{s_0}(x_{s_1})=0$. \end{itemize} Theorem \ref{main} is proved for all $s\in\mathcal{Q}$. \end{proof}

\section{Formal Markoff map} \llabel{grozarbre}

\begin{figure}[h!] \centering \input{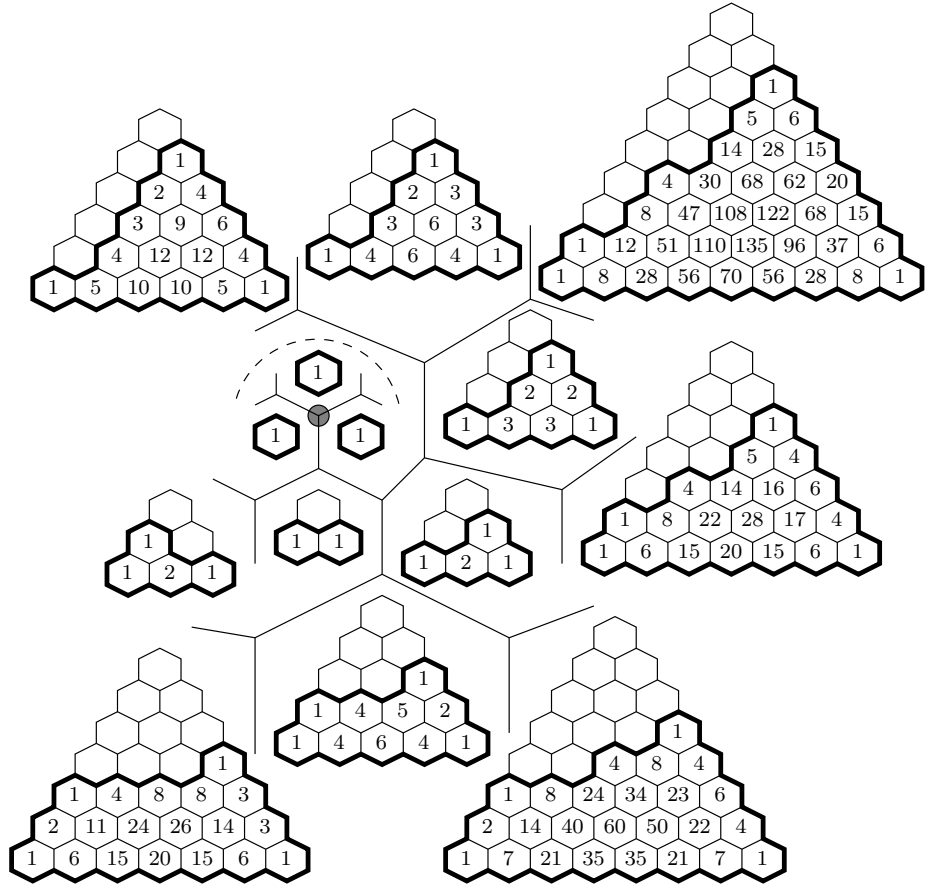}
\caption{The universal (formal) Markoff map. The integers $F_s(\cdot,\cdot)$ inside each ``bag'' add up to a Markoff number. \llabel{abeilles}} \end{figure}

Figure \ref{abeilles} shows the domains $J_s$ and the values of $F_s$ for some
of the simplest rationals $s\in\mathcal{Q}$. In each case, the points $x$ of the
affine lattice $Z_s$ have been identified with the cells of a honeycomb,
carrying the numbers $F_s(x)$. Empty cells carry $0$, by convention. Coordinates
have been tilted so that the edge $P^s_0Q^s_0$ of $J_s$ is always at the top of
$J_s$, rather than the bottom left as in Figure \ref{vertices}. The left edge of
$J_s$ consists of $p$ cells (the $P^s_i$); the right edge, of $q$ cells (the
$Q^s_j$). The single cells to the bottom left and bottom right of the ``root''
(dark spot) correspond to the exceptional cases $s=0$ and $s=\infty$. The single
cell above the root corresponds to $s=-1$; the meaning of that convention,
already apparent from the Introduction, will be re-emphasized in a moment. Observe the $1's$ in the corners of each $J_s$, just as in Lemma \ref{extremal}. It is an easy exercise (left to the reader) to prove by induction that the bottom, left, and right edges of each $J_s$ (for $s\in \mathcal{Q}-\{0,\infty\}$) always carry full lines of the Pascal triangle: if $v=(2,-2)$ then
$$F_s(P^s_i)=\binomial{p-1}{i};~F_s(Q^s_j)=\binomial{q-1}{j};
~F_s(Q^s_{q-1}+kv)=\binomial{p+q-1}{k}.$$

Notice the arrangement of the various $J_s$ in the complement $U$ of a planar
$3$-valent tree: this tree should be seen as the dual of the Farey triangulation
of $\mathbb{H}^2$, so each connected component $R_s$ of $U$ corresponds to a
horosphere centered at a rational point $s$. Each configuration $s,s_0,s_1,s'$
as in the previous section corresponds in fact to a pair of edge-adjacent
components $R_{s_0},R_{s_1}$ of $U$, together with their two common neighbors
$R_s,R_{s'}$. Since Formula (\ref{convolez}) is symmetric in $s,s'$, one may
apply it backwards to define $F_s$ for all $s$ in $\mathbb{P}^1\mathbb{Q}$
(not just $\mathcal{Q}$). This was (very) partially done in Figure
\ref{abeilles} by showing $J_{-1}=\{(-1,-1)\}$ just above the root. However, the
full picture would exhibit a $6$-fold dihedral symmetry around the root, so only
one sixth of the tree is explored to some depth in Figure \ref{abeilles}. This
$6$-fold symmetry is also the reason why honeycombs were used instead of, say,
square cells. As an exercise, the reader may prove the following formulas for
the symmetry (true for all $s\in\mathbb{P}^1\mathbb{Q}$) by induction on the tree:
$$F_{\frac{1}{s}} (\alpha,\beta)=F_{s}(\beta,\alpha)~~;~~
F_{-1-s}(\alpha,\beta)=F_{s}(-2-\alpha-\beta,\beta)$$ (The
M\"obius transformations acting on the index $s$ permute the rationals $-1,0,\infty$ while the affine transformations acting on the argument $(\alpha,\beta)$ permute the associated singletons $J_{-1},J_0, J_{\infty}$, as well as the elements of $-\Lambda$).

\section{Conjectural generalization} \llabel{generalization}
The Markoff polynomial $M=X^2+Y^2+Z^2-XYZ$ encountered in Section \ref{laurent} has degree $2$ in each variable. This is why any solution $(X,Y,Z)$ of the equation $M=0$ defines many other solutions: by considering $M$ as a polynomial of degree $2$ in, say, the variable $X$, we can always replace $X$ by the conjugate root. Thus, the free product $G$ of three copies of $\mathbb{Z}/2\mathbb{Z}$ acts naturally on the variety $M=0$ by isomorphisms. An analogous statement holds true if we replace $M$ by \emph{any} polynomial of degree $2$ in all its variable (allowing for such monomials as $X^2Y^2ZT$), and allow for actions by birational isomorphisms.

In this section, we conjecture a generalization of Theorem \ref{main} to all $N$-variable polynomials $M$ which are \emph{monic of degree $2$} in each variable. Namely, we show that certain expressions for the action of $G$ are Laurent polynomials (as in Proposition \ref{evendegrees}), and conjecture that the coefficients are positive. The coefficients of $M$ will be considered as variables themselves (noted $A_I$ below). We work over the complex field $\mathbb{C}$.

Let $N\geq 2$ be an integer, and denote by $[\![N]\!]$ the set of integers $\{1,2,\dots,N\}$. For each $I\subset [\![N]\!]$, fix a formal parameter $A_I$. Consider the Markoff-type equation in $N$ variables $X_1\dots,X_N$:
\begin{equation}\llabel{markofftype} 
\sum_{i=1}^N X_i^2+\sum_{I\subset [\![N]\!]}A_I\prod_{i\in I}X_i=0~~.
\end{equation}
Let $V\subset \mathbb{C}^N$ be the variety defined by (\ref{markofftype}). For each $k\in [\![N]\!]$ and each point $(x_1,\dots,x_N)$ of $V\cap{\mathbb{C}^*}^N$, define
\begin{equation} \begin{array}{rrl} E_k(x_1,\dots,x_N)&:=&(x_1,\dots,x_{k-1},\overline{x_k},x_{k+1},\dots,x_n) \\ \text{where }\hspace{25pt}\overline{x_k}&=&\displaystyle{\left . \left ( \sum_{i\neq k} x_i^2 + \sum_{I\subset [\![N]\!]-\{k\}}A_I\prod_{i\in I}x_i \right ) \right / x_k~~~}. \end{array} \llabel{xkbar} \end{equation}
Then $E_k$ defines a birational $\mathbb{Z}/2\mathbb{Z}$-action on $V$: indeed, $\overline{x_k}x_k$ is the product of the roots of (\ref{markofftype}), seen as a monic degree $2$ polynomial in the $k$-th variable. By letting $k$ range over $[\![N]\!]$, we obtain a birational action on $V$ by the free product $G$ of $N$ copies of $\mathbb{Z}/2\mathbb{Z}$.

Observe that the variable $A_{[\![N]\!]}$ is absent from the definition (\ref{xkbar}) of each generator $E_k$: therefore, $G$ acts on each ``level manifold'' of $\mathbb{C}^N$ defined by 
\begin{equation} B(x_1,\dots,x_N):=\left . \left ( \sum_{i=1}^N x_i^2+\sum_{I\subsetneq [\![N]\!]}A_I\prod_{i\in I}x_i \right )\right / \prod_{i=1}^N x_i~~=~~ \text{constant} \llabel{levelcurve} \end{equation}
(indeed, $B(x_1,\dots,x_N)$ is just the value of $A_{[\![N]\!]}$ for which a given point $(x_1,\dots,x_N)$ will satisfy (\ref{markofftype}), when all the $\{A_I\}_{I\subsetneq [\![N]\!]}$ are given). In particular, $B(x_1,\dots,x_N)$ is invariant under the action of $E_k$ on $\mathbb{C}^N$: therefore, the expression given in (\ref{xkbar}) for $E_k$ extends to a birational involution of $\mathbb{C}^N$ respecting $B$. Henceforward, we consider $G$ as acting on $\mathbb{C}^N$ by birational isomorphisms. 

\begin{proposition}
For each $g$ in $G$ and $x=(x_1,\dots,x_N)$ in $\mathbb{C}^N$, the coordinates of $g\cdot x$ are polynomials in the variables $\{x_i^{\pm 1}\}_{i\in [\![N]\!]}$ and $\{A_I\}_{I\subsetneq [\![N]\!]}$ with integer coefficients depending only on $g$.
\end{proposition}
\begin{remark} We conjecture that these integers are positive. Theorem \ref{main} corresponds to $N=3$ under the specialization $A_I\equiv 0$: for example, $E_1(x,y,z)=(\frac{y^2+z^2}{x},y,z)$. \end{remark}
\begin{proof}
We work by induction in $G$, using the generators $E_k$. When $g$ is the identity of $G$, we are done. Suppose the proposition is true for $g$, so that $g\cdot(x_1,\dots,x_N)=(y_1,\dots,y_N)$ where each $y_j$ is a polynomial in the $\{x_i^{\pm}\}_{i\in [\![N]\!]}$ and $\{A_I\}_{I\subsetneq [\![N]\!]}$ with integer coefficients. We must prove that the coordinates of $$E_k(y_1,\dots,y_N)=(y_1,\dots,\overline{y_k},\dots,y_N)$$ are polynomials as well, where $\overline{y_k}$ is given as in (\ref{xkbar}).
We saw that the left member $B(x_1,\dots,x_N)$ of (\ref{levelcurve}) is (formally) $E_k$-invariant for each $k\in [\![N]\!]$; therefore we must have $B(x_1,\dots,x_N)=B(y_1,\dots,y_N)$. Using (\ref{xkbar}), note that
\begin{eqnarray*}
\overline{y_k}&=&\left . \left ( \sum_{i\neq k} y_i^2 + \sum_{I\subset [\![N]\!]-\{k\}}A_I\prod_{i\in I}y_i \right ) \right / y_k \\ &=& \left . \left ( \left ( B(y_1,\dots,y_N)\prod_{i=1}^N y_i \right ) - y_k^2 - \underset{k\in I}{\sum_{I\subsetneq [\![N]\!]}} A_I \prod_{i\in I} y_i \right ) \right / y_k \\ &=& B(x_1,\dots,x_N) \left ( \prod_{i\in [\![N]\!]-\{k\}} y_i \right )-y_k - \underset{k\in I}{\sum_{I\subsetneq [\![N]\!]}} A_I \prod_{i\in I-\{k\}}y_i \end{eqnarray*}
Using the formula (\ref{levelcurve}) for $B(x_1,\dots,x_N)$,
the last expression is clearly a polynomial in the variables $\{x_i^{\pm 1}\}_{i\in [\![N]\!]}$ and $\{A_I\}_{I\subsetneq [\![N]\!]}$ with integer coefficients. This is a direct analogue of (\ref{firstdefinition}).
\end{proof}

\begin{flushright}
USC Mathematics (KAP) \\ 
3620 South Vermont Avenue \\ 
Los Angeles, CA 90089, USA \\ 
\vspace{8pt}
D\'epartement de math\'ematiques et applications, \\
\'ENS -- DMA (CNRS, UMR-8553) \\
45 rue d'Ulm, Paris Cedex 05, France \\
\texttt{Francois.Gueritaud@normalesup.org} 
\end{flushright}

\end{document}